\documentclass[12pt]{article}
\usepackage{amssymb,amsthm,amsmath,amsfonts,latexsym,tikz,hyperref,ytableau,authblk}
\usepackage[hmargin=1in,vmargin=1in]{geometry}
\usetikzlibrary{calc}

\newtheorem{thm}{Theorem}[section]
\newtheorem{prop}[thm]{Proposition}
\newtheorem{cor}[thm]{Corollary}
\newtheorem{lem}[thm]{Lemma}
\newtheorem{conj}[thm]{Conjecture}
\newtheorem{exa}[thm]{Example}

\DeclareMathOperator{\sech}{sech}
\DeclareMathOperator{\Fix}{Fix}

\newcommand{\ben}{\begin{enumerate}}
\newcommand{\een}{\end{enumerate}}
\newcommand{\ble}{\begin{lem}}
\newcommand{\ele}{\end{lem}}
\newcommand{\bth}{\begin{thm}}
\renewcommand{\eth}{\end{thm}}
\newcommand{\bpr}{\begin{prop}}
\newcommand{\epr}{\end{prop}}
\newcommand{\bco}{\begin{cor}}
\newcommand{\eco}{\end{cor}}
\newcommand{\bcon}{\begin{conj}}
\newcommand{\econ}{\end{conj}}
\newcommand{\bde}{\begin{defn}}
\newcommand{\ede}{\end{defn}}
\newcommand{\bex}{\begin{exa}}
\newcommand{\eex}{\end{exa}}
\newcommand{\barr}{\begin{array}}
\newcommand{\earr}{\end{array}}
\newcommand{\btab}{\begin{tabular}}
\newcommand{\etab}{\end{tabular}}
\newcommand{\beq}{\begin{equation}}
\newcommand{\eeq}{\end{equation}}
\newcommand{\bea}{\begin{eqnarray*}}
\newcommand{\eea}{\end{eqnarray*}}
\newcommand{\bal}{\begin{align*}}
\newcommand{\bce}{\begin{center}}
\newcommand{\ece}{\end{center}}
\newcommand{\bpi}{\begin{picture}}
\newcommand{\epi}{\end{picture}}
\newcommand{\bpp}{\begin{picture}}
\newcommand{\epp}{\end{picture}}
\newcommand{\bfi}{\begin{figure} \begin{center}}
\newcommand{\efi}{\end{center} \end{figure}}
\newcommand{\bprf}{\begin{proof}}
\newcommand{\eprf}{\end{proof}\medskip}
\newcommand{\capt}{\caption}
\newcommand{\bsl}{\begin{slide}{}}
\newcommand{\esl}{\end{slide}}
\newcommand{\bfr}{\begin{frame}}
\newcommand{\efr}{\end{frame}}

\newcommand{\comp}{\models}

\newcommand{\hqed}{\hfill \qed}

\newcommand{\eqqed}[1]{$\rule{1ex}{0ex}\hfill{\dil#1}\hfill\qed$}

\newcommand{\hs}[1]{\hspace{#1}}
\newcommand{\hso}[1]{\hspace{-1pt}}

\newcommand{\sbe}{\subseteq}

\newcommand{\setm}{\setminus}

\newcommand{\Cong}{\equiv}

\newcommand{\zh}{\hat{0}}

\newcommand{\ptn}{\vdash}

\newcommand{\case}[4]{\left\{\barr{ll}#1&\mbox{#2}\\#3&\mbox{#4}\earr\right.}

\def\<{\langle}
\def\>{\rangle}

\newcommand{\spn}[1]{\langle{#1}\rangle}

\newcommand{\ra}{\rightarrow}

\newcommand{\al}{\alpha}
\newcommand{\be}{\beta}

\newcommand{\io}{\iota}

\newcommand{\si}{\sigma}

\newcommand{\bbZ}{{\mathbb Z}}

\newcommand{\cE}{{\cal E}}

\newcommand{\cL}{{\cal L}}

\newcommand{\fS}{{\mathfrak S}}

\DeclareMathOperator{\Des}{Des}

\DeclareMathOperator{\id}{id}

\DeclareMathOperator{\Mod}{mod}

\DeclareMathOperator{\sgn}{sgn}

\newcommand{\dil}{\displaystyle}

\begin{document}
\pagestyle{plain}

\title{Generalized Euler numbers  and ordered set partitions
}
\author[1]{Bruce E. Sagan}
\affil[1]{Department of Mathematics, Michigan State University, East Lansing, MI 48824}

\date{\today\\[10pt]
	\begin{flushleft}
	                                       \\[5pt]
%https://mathscinet.ams.org/mathscinet/msc/msc2020.html
	\end{flushleft}}

\maketitle

\begin{abstract}

The Euler numbers $E_n$ have been widely studied.  A signed version of the Euler numbers of even subscript are given by the coefficients of  the exponential generating function $1/(1+x^2/2!+x^4/4!+\cdots)$.  Leeming and MacLeod introduced a generalization of the Euler numbers depending on an integer parameter $d\ge2$ where one takes the coefficients of the expansion of $1/(1+x^d/d!+x^{2d}/(2d)!+\cdots)$.  These numbers $\cE_n^{(d)}$ have been shown to have many interesting properties despite being much less studied.  And the techniques used have been mainly algebraic.  We propose a combinatorial model for the
 $\cE_n^{(d)}$ as signed sums over ordered partitions.  We show that this approach can be used to prove a number of old and new results including a recursion, integrality, and various congruences.  Our methods include sign-reversing involutions and M\"obius inversion over partially ordered sets.

\end{abstract}

\section{Introduction}
\label{int}

\begin{table}
$$
\barr{l|rrrrrrrrrr}
n 	& 0 	& 1 	& 2 	& 3 	& 4 	& 5 	&6	&7	&8	&9\\
\hline
E_n 	& 1	&1  	& 1 	& 2 	& 5	& 16 	& 61	& 272  &1385 	& 7936
\earr
$$
\capt{The Euler numbers as defined by $\tan x + \sec x$}
\label{EnTab}
\end{table}

The {\em Euler numbers}, $E_n$, can be defined  in terms of the exponential generating function
$$
\sum_{n\ge 0} E_n \frac{x^n}{n!} = \tan x + \sec x.
$$
The first few Euler number are given in Table~\ref{EnTab}.
There is a tremendous literature surrounding these constants, for example, in combinatorics and number theory.  Considering the parity of the powers of $x$ we see that
\beq
\label{sec}
\sum_{n\ge 0} E_{2n} \frac{x^{2n}}{(2n)!} =  \sec x.
\eeq

\begin{table}
$$
\barr{l|rrrrrrrrrr}
n 	& 0 	& 1 	& 2 	& 3 	& 4 	& 5 	&6	&7	&8	&9\\
\hline
\cE_n 	& 1	&0  	&- 1 	& 0 	& 5	& 0 	& -61	& 0 	 &1385 & 0
\earr
$$
\capt{The Euler numbers as defined by $2/(e^x+e^{-x})$}
\label{cEnTab}
\end{table}
The even subscripted Euler numbers have also been considered  in another context which will be amenable to generalization.
Define a sequence $\cE_n$ by
\beq
\label{cEnDef}
\sum_{n\ge0} \cE_n \frac{x^n}{n!} = \frac{2}{e^x+e^{-x}}=\frac{1}{1+x^2/2!+x^4/4!+\cdots}.
\eeq
The beginning of this sequence is displayed in Table~\ref{cEnTab}.
It is easy to see that
\beq
\label{cEE}
\cE_n =\case{(-1)^{n/2} E_n}{if $n$ is even.}{0}{if $n$ is odd.}
\eeq
Indeend, since $2/(e^x+e^{-x})$ is an even function we have the second case of~\eqref{cEE}.  As far as the first,  we can rewrite equation~\eqref{sec} as
\begin{align*}
\sum_{n\ge 0} (-1)^n E_{2n} \frac{x^{2n}}{(2n)!}
& =  \sum_{n\ge 0} E_{2n} \frac{(ix)^{2n}}{(2n)!}\\
& =  \sec( ix)\\
&=\sech(x)\\
& =\frac{2}{e^x+e^{-x}}\\
&=\sum_{n\ge 0} \cE_{2n} \frac{x^{2n}}{(2n)!}.
\end{align*}
As is done in the literature, we will also refer to the $\cE_n$ as Euler numbers and let context distinguish between them and the $E_n$.

Lehmer~\cite{leh:lrf} introduced an analogue of the Euler numbers as follows.  Let $\zeta$ be a primitive cube root of unity.  Now define the {\em Lehmer numbers}, $\cL_n$, by
\beq
\label{cWnDef}
\sum_{n\ge0} \cL_n \frac{x^n}{n!} = \frac{3}{e^x+e^{\zeta x} + e^{\zeta^2 x}}=\frac{1}{1+x^3/3!+x^6/6!+\cdots}.
\eeq
See Table~\ref{cLTab} for some specific values.
It has been shown that the $\cL_n$ also have interesting properties such as recurrences, congruences, and determinantal identities~\cite{BK:lge,KL:cpl,KP:hcn}.
Almost all of these have been derived by algebraic means such as manipulation of sums.

Looking at equations~\eqref{cEnDef} and~\eqref{cWnDef} suggests an obvious generalization.  let $d\ge 2$ be an integer and let $\zeta_d$ be a primitive $d$th root of unity.  Define the {\em generalized Euler numbers}, $\cE_n^{(d)}$, by 
\beq
\label{cEndDef}
\sum_{n\ge0} \cE_n^{(d)} \frac{x^n}{n!} = \frac{d}{e^x+e^{\zeta_d x} + e^{\zeta_d^2 x}+\cdots+e^{\zeta_d^{d-1} x}}
=\frac{1}{1+x^d/d!+x^{2d}/{2d}!+\cdots}.
\eeq
Clearly $\cE_n^{(2)} = \cE_n$ and $\cE_n^{(3)} = \cL_n$.  
These numbers were first defined by Leeming and MacLeod~\cite{LM:pge} and have since only been studied in~\cite{ges:cge,KL:cpl,LM:gen}.  Given the vast literature on Euler numbers, we feel that this generalization has been overlooked.

\begin{table}
$$
\barr{l|rrrrrrrrrr}
n 	& 0 	& 1 	& 2 	& 3 	& 4 	& 5 	&6	&7	&8	&9\\
\hline
\cL_n 	& 1	&0  	&0 	& -1	& 0	& 0 	& 19	& 0 	 &0 	 & -1513
\earr
$$
\capt{The Lehmer numbers}
\label{cLTab}
\end{table}

The purpose of the current work is to study the $\cE_n^{(d)}$ from a combinatorial viewpoint.  It is well known that the $E_n$ count alternating permutations. 
The $\cE_n^{(d)}$ have a combinatorial interpretation in terms of ordered set partitions.  Let $S$ be a set. 
Often $S$ will be the interval $[n]=\{1,2,\dots,n\}$.
 An {\em ordered set partition} of $S$ is a sequence of nonempty subsets $\pi=(B_1,B_2,\ldots,B_k)$ where $\uplus_i B_i = S$ (disjoint union).  The $B_i$ are called {\em blocks} and their order matters, while the order of the elements within each block does not. And in examples we will not write out set braces and commas unless they are needed for readability.  For example, two different ordered set partition of $[5]$ are
$$
\pi= (14,25,3) \text{ and } \si = (3,14,25).
$$
We write $\pi\comp S$ if $\pi$ is an ordered set partition of $S$.    The {\em length} of $\pi$ is the number of blocks and denoted $\ell(\pi)$.  
Both of the displayed partitions above have $\ell(\pi)=3$.  Ordered set partitions play a crucial role when considering ordered  Stirling numbers and $q$-Stirling numbers of the second kind as well as associated algebraic structures such as coinvariant algebras.  See the article of Sagan and Swanson~\cite{SS:qStB} for a history and references.

We will  write $\pi\comp_d S$ if every block $B$  of $\pi$ has $\#B$ divisible by $d$ where we use 
$\# B$ or $|B|$ to denote cardinality.  Such partitions will be called {\em $d$-divisible}.  To illustrate,
$$
\pi = (27,1346,58)\comp_2 [8].
$$
Our main tool will be the following combinatorial description of generalized  Euler numbers.
\bth
\label{cEndPtn}
For all $n\ge0$ and $d\ge2$ we have
$$
\cE_n^{(d)} = \sum_{\pi\ptn_d\hs{2pt} [n]} (-1)^{\ell(\pi)}.
$$
\eth
To illustrate, one can compute from~\eqref{cEnDef} that $\cE_4^{(2)}=\cE_4=5$.  On the other hand, the $2$-divisible partitions of $[4]$ are
$$
(1234), (12,34), (34,12), (13,24), (24,13), (14,23), (23,14).
$$
So
$$
\cE_4 = (-1)^1 + (-1)^2  + (-1)^2 + (-1)^2 + (-1)^2 + (-1)^2 + (-1)^2  = 5.
$$ 

The rest of this paper is structured as follows.
In the next section we study the  Euler numbers $\cE_{n}^{(2)}$ through the lens of ordered set partitions.  We show that they can be considered as signed sums over certain ordered partitions.  This model is then used to prove various classical results including that they are integers, alternate in sign, and satisfy a nice recursion.  Section~\ref{bpg} is devoted to showing that similar results (with similar proofs) hold for $\cE_{n}^{(d)}$ for all $d\ge2$.
In the section following that, we prove generalizations of various congruences already known for small values of $d$ either to all $d$ or to all prime $d$.
We end with some suggestions for future work.

%%%%%%%%%%%%%%%%%%%%%%%%%%%%%%%%

\section{The original Euler numbers}

We will first concentrate on the case $d=2$ where  $\cE_n^{(2)}=\cE_n$.  Several of our results and proofs will generalize easily to arbitrary $d$.  So, in those cases, we will be able to merely mention any necessary changes to get the appropriate generalization in Section~\ref{bpg}.  Our first order of business will be to prove equation~\eqref{cEndPtn} for $d=2$.  Given a power series $f(x)$ we will use the notation $[x^n/n!] f(x)$ for the coefficient of $x^n/n!$ in $f(x)$.
\bth
\label{cEnPtn}
For all $n\ge0$ we have
\beq
\label{cEnSum}
\cE_n = \sum_{\pi\comp_2\hs{2pt} [n]} (-1)^{\ell(\pi)}.
\eeq
\eth
\bprf
If $n_1+n_2+\cdots+ n_k=n$ is a composition of $n$ into positive parts then one has the corresponding multinomial coefficient
$$
\binom{n}{n_1, n_2,\ldots,n_k} = \frac{n!}{n_1! n_2!\cdots n_k!}.
$$
It is well known and easy to prove that this multinomial coefficient is the number of ordered partitions $\pi=(B_1,B_2,\ldots,B_k)\comp [n]$ with $\#B_i = n_i$ for $1\le i\le k$.  It follows that
\beq
\label{det}
\sum_{\pi\ptn_2\hs{2pt} [n]} (-1)^{\ell(\pi)} = n! \sum_{k\ge0} (-1)^k \sum_{2n_1 + 2n_2+\cdots 2n_k=n} \frac{1}{(2n_1)! (2n_2)!\cdots (2n_k)!}.
\eeq

On the other hand, recall from equation~\eqref{cEnDef} that
$$
\cE_n
= [x^n/n!] (1 + x^2/2! + x^4/4! +\cdots)^{-1}.
$$
Using the usual expansion for the inverse of a series and taking the indicated coefficient gives the same formula as for the alternating sum, so we are done.
\eprf

As an immediate corollary we get the following classical results.
\bco
For all $n\ge0$ we have 
\ben
\item[(a)] $\cE_n\in\bbZ$, and
\item[(b)] $\cE_{2n+1} = 0$.\hqed
\een
\eco

As for the $\cE_{2n}$, we can also easily obtain the following recursion.
\bpr
We have $\cE_0=1$ and for $n\ge1$
$$
\cE_{2n}=-\sum_{i=0}^{n-1} \binom{2n}{2i} \cE_{2i}
$$
\epr
\bprf
Suppose $\pi=(B_1,\ldots,B_k)$.
Consider the contribution of all such $\pi$ with $\# B_1 =2i$ to the sum~\eqref{cEnSum}.  There are $\binom{2n}{2i}$ ways to choose $B_1$.
And $\pi':=(B_2,\ldots,B_k)\comp_2 [2n]-B_1$ can be chosen in $\cE_{2n-2i}$ ways.  So the total contribution is
$-\binom{2n}{2i} \cE_{2n-2i}$ where the negative sign comes from the fact that $\pi$ has one more block than $\pi'$.  Summing over $i\in[n]$, replacing $i$ by $n-i$, and using the symmetry of the binomial coefficients finishes the proof.
\eprf

Next we would like to show that the $\cE_{2n}$ alternate in sign and count certain permutations.  Of course, this follows from equation~\eqref{cEE} and the facts that the $E_n$ are positive and enumerate alternating permutations.  But we wish to give a combinatorial proof.  We will use the standard method of sign-reversion involutions which we now briefly review.  For more information on this technique, see the text of Sagan~\cite[Section 2.2]{sag:aoc}.

\bfi
\begin{tikzpicture}
\draw (1,1) node{$(12,34)$};
\draw (1,-1) node{$(1234)$};
\draw[<->] (1,.7)--(1,-.7);
\draw (3,1) node{$(34,12)$};
\draw[->] (3.5,1.4)  arc (0:180:.5);
\draw (5,1) node{$(13,24)$};
\draw[->] (5.5,1.4)  arc (0:180:.5);
\draw (7,1) node{$(24,13)$};
\draw[->] (7.5,1.4)  arc (0:180:.5);
\draw (9,1) node{$(13,24)$};
\draw[->] (9.5,1.4)  arc (0:180:.5);
\draw (11,1) node{$(24,13)$};
\draw[->] (11.5,1.4)  arc (0:180:.5);
\draw (6,0) ellipse (7 and 3);
\draw (-1,0)--(13,0);
\draw(0,.6) node{$S^+$};
\draw(0,-.6) node{$S^-$};
\draw(-2,0) node{$S=$};
\end{tikzpicture}
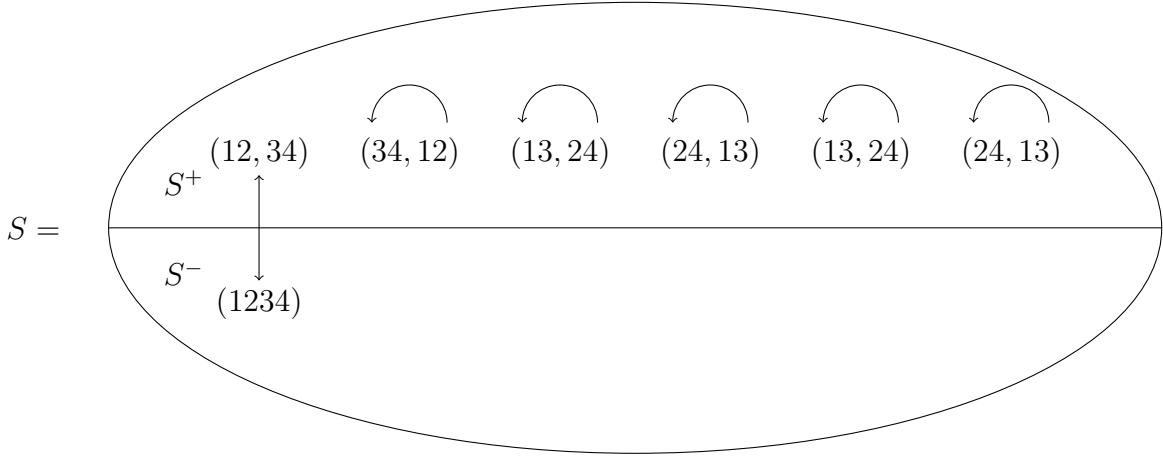
\capt{The signed set of all $2$-divisible partitions of $[4]$ and a sign-reversing involution}
 \label{ioFig}
\efi

Let $S$ be a finite set and $\io:S\ra S$ an involution, that is, $\io^2=\id$ where $\id$ is the identity map.  So $\io$ can be viewed as a permutation of $S$ whose cycle decomposition consists of $2$-cycles and fixed points.  In Figure~\ref{ioFig}, the set $S$ is all $2$-divisible partitions of $[4]$ and $\io$ is indicated by the arcs.  
So, $(12,34)$ and  $(1234)$ form a $2$-cycle and all the other ordered partitions are fixed points.
Now assume that $S$ is signed so that there is a function $\sgn:S\ra\{1,-1\}$.  We let 
$$
S^+ =\{s\in S \mid \sgn s = 1\} \text{ and } S^-=\{s\in S \mid \sgn s = -1\}.
$$
Continuing our example, we let 
$$
\sgn\pi = (-1)^{\ell(\pi)}
$$
so that $(1234)$ is the sole element with sign $-1$ and all the rest have sign $1$.  Say that $\io$ is {\em sign reversing} if, for each of its $2$-cycles $(s,t)$ we have
$$
\sgn t = -\sgn s.
$$
Our example $\io$ is clearly sign-reversing as $\sgn(12,34) = - \sgn(1234)$.
Let $\Fix\io$ be the set of fixed points of $\io$.  If $\io$ is sign reversing then  we clearly have
\beq
\label{SRI}
\sum_{s\in S} \sgn s = \sum_{s\in\Fix \io} \sgn s
\eeq
since the signs in each $2$-cycle cancel each other.  The hope is that the sum on the right will have far fewer terms and, if we are lucky, that they all have the same sign.  Figure~\ref{ioFig} illustrates the involution used in the demonstration of the next result.   A more complicated example will be found after the proof.  The use of splitting and merging to create involutions can also be used, e.g., to find cancellation-free antipodes for Hopf algebras as shown by Benedetti and Sagan~\cite{BS:ai}.

To make the connection with alternating permutations, let $\fS_n$ be the symmetric group of all permutations $\si=\si_1 \si_2\ldots \si_n$ written in $1$-line notation.  The {\em descent set} of $\si$ is
$$
\Des\si = \{i \mid \si_i>\si_{i+1}\}.
$$
Call $\si$ {\em alternating} if
$$
\Des\si =\{2i \mid 2\le 2i<n\},
$$
and  let
$$
A_n = \{\si\in\fS_n \mid \text{$\si$ is alternating}\}.
$$
It is well known that
$$
E_n = \# A_n.
$$
The next result is a restatement of~\eqref{cEE} for the even indices, but now we can give a combinatorial proof.
\bpr
\label{cEnSgn}
For all $n\ge0$ we have 
$$
\cE_{2n}=(-1)^n E_{2n}
$$
\epr
\bprf
 Consider the set
$$
\Pi_{2n}^{(2)} = \{\pi \mid \pi\comp_2 [2n]\}
$$
signed by
\beq
\label{SgnPi}
\sgn \pi =(-1)^{\ell(\pi)}.
\eeq
Appealing to Theorem~\ref{cEnPtn}, we obtain
\beq
\label{LtSum}
\sum_{\pi\in \Pi_{2n}^{(2)}} \sgn\pi = \sum_{\pi\comp_2\hs{2pt} [2n]} (-1)^{\ell(\pi)} = \cE_{2n}.
\eeq

To define the necessary involution, $\io:\Pi_{2n}^{(2)}\ra \Pi_{2n}^{(2)}$, say that a block $B_i$ of  $\pi=(B_1,\ldots,B_k)\in\Pi_{2n}^{(2)}$ is {\em splittable} if $\#B_i\ge 4$.  On the other hand, we call $B_i$ {\em mergeable} if
\ben
\item[(M1)] $\#B_i = 2$, and
\item[(M2)] $\max B_i < \min B_{i+1}$.
\een
Note the (M2) assumes that $i<\ell(\pi)$.
Note also that because of the cardinality constraints, a block can not be both splittable and mergeable.  If $\pi$ has no splittable or mergeable blocks then it is a fixed point of $\io$.  Otherwise, let $i$ be the smallest index such that the block $B_i=\{b_1<b_2<\ldots\}$ is splittable or mergeable and define
$$
\io(\pi) = 
\case{(B_1,\ldots,B_{i-1}, \{b_1,b_2\}, B_i-\{b_1,b_2\}, B_{i+1},\ldots,B_k)}{if $B_i$ is splittable,}{(B_1,\ldots,B_{i-1}, B_i \uplus B_{i+1},B_{i+2},\ldots,B_k)}{if $B_i$ is mergeable.}
$$

It is clear from~\eqref{SgnPi} that $\io$ is sign-reversing since if $\io(\pi)=\pi'$ then $\ell(\pi') = \ell(\pi)\pm 1$.  We also need to check that $\io^2=\id$.  We will show that this is the case when $\pi'$ is obtained from $\pi$ by splitting $B_i$ as the merging case is similar.  Since $b_1,b_2$ are the smallest two elements of $B_i$ we have that $\{b_1,b_2\}$ is mergeable with $B_i-\{b_1,b_2\}$.  
So we will have $\io(\pi')=\pi$ as long as the split did not create a block
$B_j$ in $\pi'$  with $j<i$  which is splittable or mergeable.  Suppose, towards a contradiction,  that such a block  did appear.  But $B_j$ could not be splittable since then $\#B_j\ge 4$ in both $\pi$ and $\pi'$ making $B_j$ splittable in $\pi$.  This contradicts the fact that $i$ was the minimal index of a splittable or mergeable block in $\pi$.  A similar argument shows that $B_j$ could not be mergeable because it would have also  been mergeable in $\pi$.

Now suppose $\pi\in\Fix\io$.  Then $\pi$ can not have any splittable blocks which means $\# B_i=2$ for all blocks $B_i$ of $\pi$.  Since $\pi\comp_2 [2n]$ this implies that $\ell(\pi) = n$  and so 
$\sgn\pi = (-1)^n$.  Thus
$$
\sum_{\pi\in\Fix\io} \sgn\pi = \sum_{\pi\in\Fix\io} (-1)^n = (-1)^n \#\Fix\io.
$$
Comparing this expression with~\eqref{LtSum} and using~\eqref{SRI} gives
$$
\cE_{2n} = (-1)^n \#\Fix\io.
$$
So, to finish the prove, we just need to show that 
\beq
\label{FixA}
\#\Fix\io = \# A_{2n}.
\eeq

We have already seen that if $\pi=(B_1,\ldots,B_n)\in\Fix\io$ then  $\#B_i = 2$ for all $i$ which makes no $B_i$ splittable. To make sure it is not mergeable we can not violate (M1).  So (M2) must be false for every pair of adjacent blocks.  Now map $\Fix\io\ra A_{2n}$ by sending $\pi$ to the permutation $\si=\si_1\ldots \si_n$ obtained  by writing each $B_i$ in increasing order and concatenating the resulting $2$-element permutations.  Condition (M2) being false for $\pi$ is equivalent to $\si$ being alternating.  So this map is a bijection which completes the proof of~\eqref{FixA} and of the proposition.
\eprf

To illustrate the involution of the proof, suppose 
$$
\pi = 2,9/4,11/1,3,5,6/7,8/10,12.
$$
Now $\{2,9\}$ has too few elements to be splittable.  And it is not mergeable with $\{4,11\}$ since
$$
\max\{2,9\}= 9> 4 =\min\{4,11\}.
$$
Similarly, $\{4,11\}$ is neither splittable nor mergeable.  But $\{1,3,5,6\}$ is splittable since it has (at least) $4$ elements.  Thus
$$
\pi'=\io(\pi) = 2,9/4,11/1,3/5,6/7,8/10,12.
$$
Note that the fact that $\{7,8\}$ is mergeable with $\{10,12\}$ in $\pi$ is irrelevant since $\{1,3,5,6\}$ comes earlier in the partition.  To compute $\io(\pi')$, we see that the first two blocks are neither splittable nor mergeable for the same reasons as in $\pi$.  But $\{1,3\}$ can be merged with $\{5,6\}$ so that $\io(\pi')=\pi$ as desired.

%%%%%%%%%%%%%%%%%%%%%%%%%%%%%%%%

\section{Basic properties of generalized  Euler numbers}
\label{bpg}

In this section we will study the $\cE_n^{(d)}$ for general $d$.  We start by recording analogues of the results from the previous section.  We also need the definition that a permutation $\si\in\fS_n$ is
{\em $d$-alternating} if
$$
\Des\si=\{di \mid d\le di <n\}.
$$
We also let
$$
A_n^{(d)} = \{\si\in\fS_n \mid \text{$\si$ is $d$-alternating}\}.
$$
\bth
\label{BigThm}
For all $d\ge2$ we have the following.
\ben
\item[(a)]  For all $n\ge0$ we have
\beq
\label{cEndPtn}
\cE_n^{(d)} = \sum_{\pi\comp_d\hs{2pt} [n]} (-1)^{\ell(\pi)}.
\eeq
\item[(b)] For all $n\ge0$ we have $\cE_n^{(d)} \in\bbZ$.
\item[(c)] For all $n\ge0$ we have $\cE_n^{(d)} =0$ if $n$ is not a multiple of $d$.
\item[(d)] We have $\cE_0^{(d)} =1$ and for $n\ge1$
$$
\cE_{dn}^{(d)} =-\sum_{i=0}^{n-1}\binom{dn}{di} \cE_{di}^{(d)}.
$$
\item[(e)]  For all $n\ge0$ we have 
$$
\cE_{dn}^{(d)} = (-1)^n\# A_{dn}^{(d)}.
$$
\een
\eth
\bprf
In all cases, the proofs of these statements are simple modifications of the demonstrations when $d=2$.  So, we will just indicate how these changes are applied to obtain (e).

The set for the involution is
$$
\Pi_{dn}^{(d)} = \{ \pi \mid \pi \comp_d [dn]\}.
$$
And the sign is exactly the same as in definition~\eqref{SgnPi}.  As far as the involution $\io:\Pi_{dn}^{(d)}\ra\Pi_{dn}^{(d)}$ itself, we call block $B_i$ 
of $\pi=(B_1,\ldots,B_k)$
{\em splittable} if $\# B_i\ge 2d$ or {\em mergeable} if
\ben
\item[(M1')] $\#B_i = d$, and
\item[(M2')] $\max B_i < \min B_{i+1}$.
\een
Now $\io$ is defined by splitting off the smallest $d$ elements of a splittable block or taking the disjoint union of a mergeable block with the following block, whichever comes first.  The partition $\pi$ is left fixed if no such block exists.  The reader should now be able to fill in the rest of the details.
\eprf

%%%%%%%%%%%%%%%%%%%%%%%%%%%%%%%%

\section{Congruences for generalized Euler numbers}

We will now derive some congruences for generalized Euler numbers.  
Proofs of similar results in the literature are algebraic while ours are combinatorial.

Our first theeorem  contains a result of  Leeming and MacLeod modulo $2$ as part (a).  But our technique works for any modulus, although the expressions become increasingly more complicated.  To illustrate the method, we have provided a full demonstration for mod $3$ in part (b).
\bth
Suppose $d\ge2$ and $n\ge0$ are arbitrary.
\ben
\item[(a)]  We have
\beq
\label{cEndMod2}
\cE_{dn}^{(d)} \Cong 1\  (\Mod 2).
\eeq
\item[(b)]  We have
$$
\cE_{dn}^{(d)} \Cong -1 + \sum_{k=1}^{n-1} \binom{dn}{dk}\ (\Mod 3).
$$
\een
\eth
\bprf
For (b), consider the action of $C_3$ on a $\pi\in\Pi_{dn}^{(d)}$ which fixes the partitions with at most $2$ blocks. And if $\pi=(B_1,B_2,B_3,B_4,\ldots,B_k)$ with $k\ge 3$ then
$$
(1,2,3)\pi = (B_2,B_3,B_1,B_4,\ldots,B_k).
$$
It follows that for $\pi$ with at least $3$ blocks we have $|C_3\pi|=3$.
 Note also that for any $g\in C_3$ and any $\pi\in\Pi_{pn}^{(p)}$ we have $\ell(\pi)=\ell(g\pi)$.  So, if $\ell(\pi)\ge3$ then
$$
\sum_{g\in C_3} (-1)^{\ell(g\pi)} = 3\cdot (-1)^{\ell(\pi)} \Cong 0\ (\Mod 3).
$$
Thus, appealing to equation~\eqref{cEndPtn},
\begin{align*}
\cE_{dn}^{(d)} 
&= \sum_{\pi\comp_d\hs{2pt} [dn]} (-1)^{\ell(\pi)}\\
& \Cong (-1)^{\ell([dn])} + \sum_{\pi\comp_d\hs{2pt} [dn] \atop \ell(\pi)=2} (-1)^2\ (\Mod 3)\\
& = -1 + \sum_{k=1}^{n-1} \binom{dn}{dk}
\end{align*}
as desired.
\eprf

Our next theorem is a generalization to an arbitrary prime $p$ of a result which was known
for $p=2$ and $3$.  The former follows from a congruence of Stern~\cite{ste:tez}.  The latter was demonstrated in a recent paper of Komatsu and Liu~\cite{KL:cpl}, although results for larger modulus had already been proved in the original paper of Leeming and MacLeod~\cite{LM:pge}. 
Our proof will use the technique of M\"obius inversion over a partially ordered set (poset) which we will now review.  More information about this method can be found in the texts of Sagan~\cite{sag:aoc} or Stanley~\cite{sta:ec1}.

Let $P$ be a poset with a unique minimal element $\zh$.  The {\em M\"obius function} of $P$ is the map $\mu:P\ra\bbZ$ defined recursively by
$$
\mu(x) = \case{1}{if  $x=\zh$,}{\dil-\sum_{y<x}\mu(y)}{otherwise.}
$$
If $P$ is the lattice of divisors of a positive integer $n$ then this function reduces to the unsual M\"obius function from number theory.  Its importance is that one can use this function to invert sums.
\bth[M\"obius Inversion Theorey]
Let $P$ be a finite poset with a $\zh$, $V$ a real vector space, and $\al,\be:P\ra V$ two functions. Then
\beq
\label{MuHyp}
\al(x) =\sum_{y\ge x} \be(y) \text{ for all $x\in P$} 
\eeq
implies that

\medskip

\eqqed{ \be(\zh) = \sum_{y\in P} \mu(y)\al(y).}
\eth

The proof of our next result  will combine group actions and M\"obius inversion, a method which can be used to prove many congruences, see~\cite{sag:cag}.

\bfi
\begin{tikzpicture}
\draw(0,0) node{$\spn{e}$};
\draw(-5,2) node{$\spn{g}$};
\draw(-3,2) node{$\spn{gh}$};
\draw(-1,2) node{$\spn{gh^2}$};
\draw(1,2) node{$\cdots$};
\draw(3,2) node{$\spn{gh^{p-1}}$};
\draw(5,2) node{$\spn{h}$};
\draw(0,4) node{$\spn{g,h}$};
\draw (-.5,.2)--(-4.5,1.8)  (-.3,.2)--(-2.5,1.8) (-.1,.2)--(-1,1.8)
(.5,.2)--(4.5,1.8)  (.3,.2)--(2.4,1.7) 
(-4.5,2.2)--(-.5,3.8) (-2.5,2.2)--(-.3,3.8) (-1,2.2)--(-.1,3.8)
(4.5,2.2)--(.5,3.8) (2.5,2.2)--(.3,3.8)
;
\end{tikzpicture}
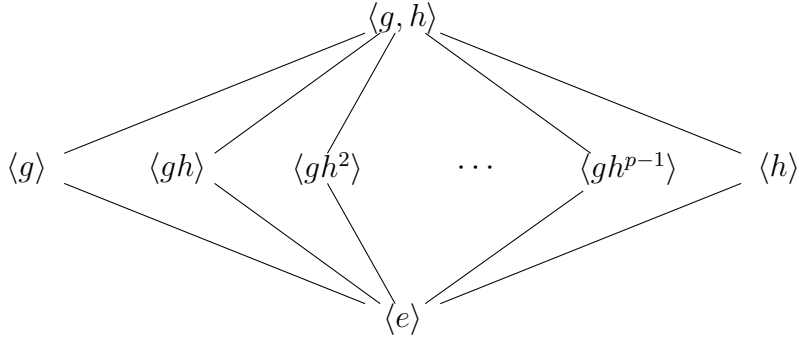
\capt{The subgroup lattice $\cL$ for $C_p\times C_p$ \label{lattice}}
\efi
\bth
\label{p^2Thm}
For  $p$ a prime we have
\beq
\label{p^2Eq}
\cE_{pn}^{(p)} \Cong (-1)^n\  (\Mod p^2).
\eeq
\eth
\bprf
Our proof will proceed in three stages.  We will first use a group action and M\"obius inversion to obtain a congruence~\eqref{HcL} for certain signed sums over stabilizers of ordered partitions.  We will then use a sign-reversing involution to  reduce the number of terms which need to be considered in the sums.  The fixed points will be in bijection with the elements of  
$\Pi_{p(n-1)}^{(p)}$ and of $\Pi_{p(n-2)}^{(p)}$  which will permit us to use induction on $n$.  Note that the result of the theorem is trivial in the base cases of  $n=0$ or $1$.

Consider the cycles $g=(1,2,\ldots,p)$ and $h=(p+1,p+2,\ldots,2p)$ in the symmetric group $\fS_{pn}$ where $n\ge2$.  The group we will use is the product 
$G=C_p \times C_p=\spn{g} \times\spn{h}$  where the angle brackets denote the group generated by an element.  
This group acts on $\Pi_{pn}^{(p)}$ by permuting the elements of the blocks according to the cycles $g$ and $h$ and their powers.
If $\pi$ is an ordered set partition then we let $G_\pi$ be the stabilizer of $\pi$.
Let $\cL$ be the lattice of subgroups of $G$ ordered by containment, see
    Figure~\ref{lattice} where $e$ is the identity element.  From the diagram it is clear that for a subgroup $H\le G$ we have
\beq
\label{muH}
\mu(H) = \begin{cases}
1  	&\text{if $H=\spn{e}$,}\\
-1  	&\text{if $H=\spn{g}, \spn{gh},\ldots, \spn{h}$,}\\
p 	&\text{if $H=\spn{g,h}$.}\\
\end{cases}
\eeq

We now define the desired functions.  Given a set of of ordered partitions $\Pi$ we let
$$
S(\Pi) = \sum_{\pi\in\Pi} (-1)^{\ell(\pi)}.
$$
Note that by equation~\eqref{cEndPtn}, we have
\beq
\label{SE}
S(\Pi_{dn}^{(d)}) = \cE_{dn}^{(d)}.
\eeq
Finally, let $\al,\be:\cL\ra\bbZ$ we given by
\beq
\label{alH}
\al(H) = S(\pi \mid G_\pi\ge H)
\eeq
and
$$
\be(H) = S(\pi \mid G_\pi= H).
$$
It is clear from the definitions that~\eqref{MuHyp} is satisfied.  So  the conclusion of the M\"obius Inversion Theorem holds and we will compute each term.

As far as $\be(\zh)$, if $H_\pi=e$ then $|G\pi|=p^2$.  Furthermore, all elements of $G\pi$ have the same number of blocks and therefore all contribute 
$(-1)^{\ell(\pi)}$ to the sum.  But this means the the total contribution of $G\pi$ is zero modulo $p^2$.  So, by M\"obius inversion,
\beq
\label{HcL}
0 \Cong \sum_{H\in\cL} \mu(H) \al(H)\ (\Mod p^2).
\eeq

To compute $\al(e)$, note that the stabilizer of every $\pi\in\Pi_{pn}^{(p)}$ contains $e$.  Appealing to~\eqref{SE} gives
\beq
\label{al(e)}
\al(e) =S(\pi\in\Pi_{pn}^{(p)}) = \cE_{pn}^{(p)}.
\eeq

Rather than compute the rest of the $\al(H)$ directly, we will employ a sign-reversing involution to cancel many of the terms.  Let us consider $\al(\spn{g})$.
Note that $G_\pi\geq \spn{g}$ if and only if the elements of $[p]$ are all in the same block of $\pi$.  Suppose $\pi=(B_1,\ldots, B_k)$ and that $[p]\sbe B_i$.
The involution $\io$ is defined by
$$
\io(\pi) =
\begin{cases}
(B_1,\ldots,B_{i-1},[p],B_i\setm [p],B_{i+1},\ldots,B_k) &\text{if $B_i\supset [p]$},\\
(B_1,\ldots,B_{i-1},B_i\uplus B_{i+1},B_{i+2},\ldots,B_k &\text{if $B_i=[p]$ for $i<k$},\\
(B_1,B_2,\ldots,B_k) &\text{if $B_k=[p]$}.
\end{cases}
$$
In other words, if the block $B_i$ containing $[p]$ contains other elements, then $[p]$ is split off as its own block and placed directly before what remains of $B_i$.  If $B_i=[p]$ but is not the last block of $\pi$ then it merges with the block after it.  And if the last block is $[p]$ then $\pi$ is a fixed point.  
Now the split and merge options are inverses of each other and they change $\ell(\pi)$ by exactly one.  So $(-1)^{\ell(\pi)} + (-1)^{\ell(\io(\pi))} = 0$ and such pairs can be ignored.  Furthermore, removing the final $[p]$ gives a bijection $\pi\leftrightarrow\pi'$  between the fixed points of $\io$ and the elements of 
$\Pi_{p(n-1)}^{(p)}$ where $\ell(\pi') = \ell(\pi)-1$.  So, by equation~\eqref{cEndPtn}, the final sum is
\beq
\label{al(g)}
\al(\spn{g}) = -\sum_{\pi'\in\Pi_{p(n-1)}^{(p)}} (-1)^{\ell(\pi')}\  = -\cE_{p(n-1)}^{(p)}.
\eeq
Clearly $\al(\spn{h})$ is given by the same expression.

It is easy to see that the subgroups $\spn{gh^i}$ for $i\in[p-1]$ as well as $\spn{g,h}$ all stabilize the same set of partitions, namely those where the elements 
of $[p]$ are in a single block and $p+1,p+2,\ldots,2p$ are in a single block (not necessarily the same as the block for $[p]$).  Using arguments similar to those in the previous paragraph and induction we obtain, for any of these subgroups $H$,
\beq
\label{al(H)}
\al(H) = \cE_{p(n-2)}^{(p)}
\eeq
Plugging~\eqref{muH}, \eqref{al(e)}, \eqref{al(g)}, and~\eqref{al(H)} into equation~\eqref{HcL} and using induction on $n$ we obtain
\begin{align*}
0
&\Cong\cE_{pn}^{(p)} - 2( -\cE_{p(n-1)}^{(p)})-(p-1)\cE_{p(n-2)}^{(p)} + p\ \cE_{p(n-2)}^{(p)}\  (\Mod p^2)\\
&\Cong \cE_{pn}^{(p)} + 2(-1)^{n-1} +(-1)^{n-2}\ (\Mod p^2)
\end{align*}
Solving for $\cE_{pn}^{(p)}$ completes the proof.
\eprf

Sometimes we can improve the modulus in equation~\eqref{p^2Eq}.  The special case when $p=3$ in the following theorem was proven by Leeming and MacLeod~\cite{LM:pge}.  It follows immediately from equations~\eqref{cEndMod2} and~\eqref{p^2Eq} as well as the fact that $2$ is relatively prime to any prime $p\ge3$.
\bth
For $p\ge3$ a prime we have

\medskip

\eqqed{
\cE_{pn}^{(p)} \Cong (-1)^n\  (\Mod 2p^2).
}
\eth

%%%%%%%%%%%%%%%%%%%%%%%%%%%%%%%%%%%%%%%%%%

\section{Future work}

The study of the generalized Euler numbers is still in its infancy.  So, there is much more to do.  We mention three avenues for future research here.

\ben
\item  Much of the work on generalized Euler numbers has concentrated on algebraic proofs of congruences~\cite{ges:cge,LM:pge,LM:gen}.  It would be interesting to see how many of them can be proven using methods such as sign-reversing involutions, M\"obius inversion, and other combinatorial techniques.
\item  There are close connections between the original Euler numbers and Bernoulli numbers.  How do these carry over to the generalized case?
\item  There are nice determinantal expressions for the original Euler numbers and the Lehmer numbers in the literature.  In fact, using an analogue of equation~\eqref{det} it can be proven that
$$
\cE_{dn}^{(d)} = (-1)^n (dn)!
\left|
\frac{1}{((i-j+1)d)!}
\right|_{1\le i,j\le n}
$$
where a term in the determinant is considered to be zero if $i-j+1<0$.  It would be very interesting to give a combinatorial proof of this result, perhaps by using the Lindstr\"om-Gessel-Viennot lattice path method~\cite{lin:vri,GV:bdp}.
\een

%\bibliographystyle{plain}

%\nocite{*}
%\bibliographystyle{abbrvnat}

\nocite{*}
\bibliographystyle{alpha}

\begin{thebibliography}{Leh35}

\bibitem[BK19]{BK:lge}
Rupam Barman and Takao Komatsu.
\newblock Lehmer's generalized {E}uler numbers in hypergeometric functions.
\newblock {\em J. Korean Math. Soc.}, 56(2):485--505, 2019.

\bibitem[BS17]{BS:ai}
Carolina Benedetti and Bruce~E. Sagan.
\newblock Antipodes and involutions.
\newblock {\em J. Combin. Theory Ser. A}, 148:275--315, 2017.

\bibitem[Ges83]{ges:cge}
Ira~M. Gessel.
\newblock Some congruences for generalized {E}uler numbers.
\newblock {\em Canad. J. Math.}, 35(4):687--709, 1983.

\bibitem[GV85]{GV:bdp}
Ira Gessel and G\'erard Viennot.
\newblock Binomial determinants, paths, and hook length formulae.
\newblock {\em Adv. in Math.}, 58(3):300--321, 1985.

\bibitem[KL25]{KL:cpl}
Takao Komatsu and Guo-Dong Liu.
\newblock Congruence properties of lehmer-euler numbers, 2025.
\newblock Preprint {\texttt{arXiv:2501.01178}}.

\bibitem[KP21]{KP:hcn}
Takao Komatsu and Ram~Krishna Pandey.
\newblock On hypergeometric {C}auchy numbers of higher grade.
\newblock {\em AIMS Math.}, 6(7):6630--6646, 2021.

\bibitem[Leh35]{leh:lrf}
D.~H. Lehmer.
\newblock Lacunary recurrence formulas for the numbers of {B}ernoulli and
  {E}uler.
\newblock {\em Ann. of Math. (2)}, 36(3):637--649, 1935.

\bibitem[Lin73]{lin:vri}
Bernt Lindstr\"om.
\newblock On the vector representations of induced matroids.
\newblock {\em Bull. London Math. Soc.}, 5:85--90, 1973.

\bibitem[LM81]{LM:pge}
D.~J. Leeming and R.~A. MacLeod.
\newblock Some properties of generalized {E}uler numbers.
\newblock {\em Canadian J. Math.}, 33(3):606--617, 1981.

\bibitem[LM83]{LM:gen}
D.~J. Leeming and R.~A. MacLeod.
\newblock Generalized {E}uler number sequences: asymptotic estimates and
  congruences.
\newblock {\em Canad. J. Math.}, 35(3):526--546, 1983.

\bibitem[Sag85]{sag:cag}
Bruce~E. Sagan.
\newblock Congruences via abelian groups.
\newblock {\em J. Number Theory}, 20(2):210--237, 1985.

\bibitem[Sag20]{sag:aoc}
Bruce~E. Sagan.
\newblock {\em Combinatorics: the art of counting}, volume 210 of {\em Graduate
  Studies in Mathematics}.
\newblock American Mathematical Society, Providence, RI, 2020.

\bibitem[SS24]{SS:qStB}
Bruce~E. Sagan and Joshua~P. Swanson.
\newblock {$q$}-{S}tirling numbers in type {$B$}.
\newblock {\em European J. Combin.}, 118:Paper No. 103899, 35, 2024.

\bibitem[Sta12]{sta:ec1}
Richard~P. Stanley.
\newblock {\em Enumerative combinatorics. {V}olume 1}, volume~49 of {\em
  Cambridge Studies in Advanced Mathematics}.
\newblock Cambridge University Press, Cambridge, second edition, 2012.

\bibitem[Ste75]{ste:tez}
M.~Stern.
\newblock Zur {T}heorie der {E}ulerschen {Z}ahlen.
\newblock {\em J. Reine Angew. Math.}, 79:67--98, 1875.

\end{thebibliography}

\end{document}